\documentclass{article}
\usepackage[top=30truemm,bottom=30truemm,left=25truemm,right=25truemm]{geometry}
\usepackage{amsmath,amssymb,amsfonts,mathtools,accents,bm,mathrsfs,bibentry,lscape,multirow}
\usepackage{enumerate}
\usepackage{tabu}
\newtheorem{theorem}{Theorem}

\newtheorem{cor}{Corollary}

\newtheorem{proposition}{Proposition}

\title{On tail inference in iid settings with nonnegative extreme value index}
\author{
Taku MORIYAMA\\
School of Data Science, Yokohama City University}
\date{}

\begin{document}
\maketitle

\begin{abstract}
In extreme value inference it is a fundamental problem how the target value is required to be extreme by the extreme value theory. In iid settings this study both theoretically and numerically compares tail estimators, which are based on either or both of the extreme value theory and the nonparametric smoothing. This study considers tail probability estimation and mean excess function estimation. 

This study assumes that the extreme value index of the underlying distribution is nonnegative. Specifically, the Hall class or the Weibull class of distributions is supposed in order to obtain the convergence rates of the estimators. This study investigates the nonparametric kernel type estimators, the fitting estimators to the generalized Pareto distribution and the plug-in estimators of the Hall distribution, which was proposed by Hall and Weissman (1997). In simulation studies the mean squared errors of the estimators in some finite sample cases are compared. 
\end{abstract}

\section{Introduction}

The occurrence probabilities of extremely rare events never happened before is known to be not exactly zero in many cases, but it is estimated to be exactly zero without any prior information. That means fully nonparametric approaches basically cannot be applied to such a extreme value inference. 

If a few assumptions are supposed on the tail of the underlying distribution $F$, extreme value theory provides ways for the inference. Specifically, the excess distribution $F_u$ for a threshold $u$ can be  approximated by the generalized Pareto distribution (GPD) $H$ if $F$ belongs to the maximum domain of attraction. It is also known to be the Pickands-Balkema-de Haan theorem, and most distributions belong to the domain. Utilizing the approximation for tail probability estimation, it is estimated by the formula $F(x) = F(u) + \{1-F(u)\} F_u(x-u)$, which is considered as the piecing-together approach (Smith 1987; Maiboroda and Markovich 2004; Macdonald et al. 2011). If the extreme value index is restricted to being positive, the asymptotic tail probability has a general form known as the Hall class of distributions introducing later. By applying plug-in principle to the estimation of the unknown parameters, the tail probability is estimated (Hall and Weissman 1997). This is the so-called plug-in approach. 

In short, there are mainly three ways in tail probability estimation: nonparametric, piecing-together and plug-in approaches. This study compares three type of estimators and concerns with a comparison between extreme value based approach (parametric) and nonparametric. In particular, how 'extreme' the target value is required by the extreme value theory is a fundamental question (see e.g. Smith 1987; Hall and Weissman 1997; Moriyama 2024; Moriyama 2024). The aim of this study is to obtain asymptotic properties of the estimators and compare the accuracy in iid settings. Specifically, it is supposed that $F$ belongs to either the Hall class or the Weibull class. 

Let $X_1,X_2,\cdots, X_n$ be independent and identically distributed random variables with a continuous distribution function $F$, where $n$ is supposed to be sufficiently large. This study considers the pointwise estimate of the tail probability $\overline{F}(x) := 1 -F(x)$ and investigates asymptotic properties of the tail probability estimators, where $x := x_n \to \infty$ as $n \to \infty$ is supposed. Each asymptotic mean squared error (MSE) is given in Section 2. A comparative study is given in Section 3. Next, this study investigates mean excess function (MEF) estimation defined as
$$e(u) := \int_{0}^{\infty} x f_u(x) {\rm d}x, ~~~\text{where} ~~~ f_u(x):= \frac{f(x+u)}{1-F(u)}.$$
MEF is also called mean residual life and is the function of the expected value of the excess distribution. MEF exists only if the extreme value index satisfies $\gamma<1$. Basic properties are provided in Embrechts et al. (1997) and the limiting behavior is investigated by Bradley and Gupta (2003). 

Necir et al. (2010) studied the related risk indicator, conditional tail expectation for the heavy tailed data. Excess distribution estimators are compared in Moriyama (2024). Asymptotic properties of the MEF estimators are provided in Section 4 and Section 5 discusses the comparative study on the estimators. The proof of Theorem 1 is given in the Appendix.

\section{Tail probability estimation}
\subsection{The kernel type estimator as a nonparametric approach}

The kernel type estimator $\widehat{F}(x)$ is a conventional method for probability estimation, where nonparametric density estimators are compared in heavy tailed cases(Takada 2008). $\widehat{F}(x)$ is given by the sample mean of $\{W(h^{-1} (x- X_i))\}_{i=1}^n$, where $W$ is the cumulative distribution of a symmetric density function $w$ called the kernel function. $h$ needs to satisfy $h \to 0$. Then, the MSE is given by the following corollary, which is a consequence of Theorem 1 in Moriyama (2024).

\begin{cor}
Suppose $F$ is continuously twice differentiable at $x$ and $\int z^2 w(z) {\rm d}z <\infty$. Then,
\begin{align*}
\mathbb{E}[(F(x) -\widehat{F}(x))^2] \sim \Biggm\{h^{2} \frac{f'(x)}{2} \int z^2 w(z) {\rm d}z  \Biggm\}^2 + \frac{\overline{F}(x)}{n}.
\end{align*}
If $\vert\int z W(z) w(z) {\rm d}z\vert < \infty$ and $n^{-1}\{f'(x)\}^{-2} f(x) \to 0$, the asymptotically optimal bandwidth is
$$h= n^{-1/3} \left(\frac{2f(x) \int z W(z) w(z) {\rm d}z}{\left(f'(x) \int z^2 w(z) {\rm d}z\right)^2}\right)^{1/3}.$$
\end{cor}
The piecing-together approach introduced next assumes that $F$ belongs to the maximum domain of attractions. Moreover, the plug-in approach needs to specify the parametric class. On the other hand, the nonparametric approach does not need such assumption.

The optimal bandwidth in the sense of a global error (e.g. mean integrated squared error) has been researched well. Examples of the estimators include the plug-in type (Altman and L\'{e}ger 1995) and the cross-validation-based type (Bowman et al. 1998). However, it is obvious that the estimators are not suitable for the tail probability estimation. At least not the mean integrated squared error but the mean weighted integrated squared error (see e.g. Altman and L\'{e}ger 1995), which puts on weight the distribution tail. If we seek the pointwise optimal bandwidth, it is necessary to estimate the extreme value index, which affects the convergence rate of the kernel type estimator. 

Let us introduce {\rm (i)} the Hall class of distributions (Hall and Welsh 1984), which satisfies $\alpha>0, \beta\ge2^{-1}, A>0, B\neq0$ and 
$$x^{\alpha+\beta} \{1 - F(x) - A x^{-\alpha} (1 + B x^{-\beta} ) \} \to 0~~~ {\rm as}  ~~~ x \to \infty.$$
Then, the asymptotically (pointwise) optimal bandwidth is
$$h_{A,\alpha} := n^{-1/3} \{2 A^{-1} \alpha^{-1}(\alpha+1)^{-2} x^{\alpha+3}\}^{1/3} \left(\frac{\int z W(z) w(z) {\rm d}z}{\left(\int z^2 w(z) {\rm d}z\right)^2}\right)^{1/3}.$$
Estimating the parameters $(A,\alpha)$ by the $r$ largest order statistics as follows
$$\widehat{\alpha} := \left\{r^{-1} \sum_{j=1}^r \ln X_{n,n-j+1} - \ln X_{n, n-r} \right\}^{-1}, ~~~ \widehat{A} := n^{-1}r (X_{n, n-r})^{1/\widehat{\gamma}}$$
(see Hall 1982), the bandwidth estimator is given by $h_{\widehat{A},\widehat{\alpha}}$, where $r$ is supposed to be of order  $n^{2\beta/(2\beta+\alpha)}$.

Since $(\widehat{\alpha} - \alpha) = O_P(n^{-\beta/(2\beta+\alpha)})$ under the assumption, suppose $(\ln x) (\widehat{\alpha} - \alpha) = o_P(1)$ i.e. $x=o(\exp(n^{\beta/(2\beta+\alpha)}))$. Then, it holds that 
$$\frac{x^{\widehat{\alpha}}}{x^{\alpha}}= 1 + (\ln x) (\widehat{\alpha} - \alpha) + (\ln x) o_P(\widehat{\alpha} - \alpha).$$ 
Combining $(\widehat{A} - A) = O_P(n^{-\beta/(2\beta+\alpha)} \ln n)$, we have the following result, which is the consequence of Hall (1982).
\begin{cor}
Suppose $F$ belongs to {\rm (i)} the Hall class, $r =O(n^{2\beta/(2\beta+\alpha)})$ and $x=o(\exp(n^{\beta/(2\beta+\alpha)}))$. Then, under the assumptions of Corollary 1, 
$$h_{A,\alpha}^{-1} (h_{\widehat{A},\widehat{\alpha}} -h_{A,\alpha}) = O_P(n^{-\beta/(2\beta+\alpha)} (\ln x + \ln n) ).$$
\end{cor}
However, we note that the pointwise bandwidth makes the estimated tail probability function lose its monotonicity. 

\subsection{The piecing-together approach}

The piecing-together is based on the Pickands-Balkema-De Haan theorem, which states that the excess distribution  (ED) $F_u(x)$ converges to the generalized Pareto distribution (GPD) $H_{\gamma, c}(x)$ as $u\uparrow \sup({\rm supp}(f))$, where 
\begin{align*}
H_{\gamma, c}(x) :=& H_{\gamma}(c^{-1}{x}) ~~~ {\rm for} ~~ 1+{\gamma} c^{-1}{x} > 0, \\
H_{\gamma}(x) :=&1 -
\begin{dcases}
(1+\gamma x)^{-1/\gamma} ~~~ &1+\gamma x>0 ~~~{\rm and} ~~~ \gamma \in \mathbb{R} \setminus \{0\} \\
\exp(-x) ~~~ &x \in \mathbb{R} ~~~{\rm and} ~~~ \gamma =0.
\end{dcases} 
\end{align*}
If $F$ belongs to {\rm (i)} the Hall class or {\rm (ii)} the following Weibull class of distributions, which satisfies
{\rm (ii)} $\kappa>0, C>0$ and 
$$\exp(Cx^{\kappa})\{1 - F(x) - \exp(-Cx^{\kappa})\} \to 0~~~ {\rm as}  ~~~ x \to \infty,$$
respectively. Then, the limiting GPD of ED is a Fr\'{e}chet or Gumbel type, and so $F_u(x)$ is approximated by $H_{\bm{\gamma}}(x)$, where $\bm{\gamma}:=(\gamma, c_u)$
\begin{align*}
\gamma :=
\begin{dcases}
\alpha^{-1} ~~~ &{\rm for ~ the ~ Hall ~ class} \\
0 ~~~ &{\rm for ~ the ~ Weibull ~ class}.
\end{dcases} ~~~~~~
c_u :=
\begin{dcases}
\gamma u ~~~ &{\rm for ~ the ~ Hall ~ class} \\
C^{-1} \kappa^{-1} u^{1-\kappa} ~~~ &{\rm for ~ the ~ Weibull ~ class}.
\end{dcases}
\end{align*}
under $\kappa\le 1$. 

The tail probability can be estimated by the following approximation
$$F(x) \sim F(u) + \{1-F(u)\} H_{\bm{\gamma}}(x-u).$$
This study supposes $u_0=F^{-1}(1-n^{-1}N_0)$. Then, the tail estimator constructed based on the extreme value theory has the following form: 
$$\widehat{H}(x) := (1-n^{-1}N) + \{n^{-1}N\} H_{\widehat{\bm{\gamma}}}(x- u),$$
where $u$ and $N$ are random variables. The unknown parameters $(\gamma, c_u)$ are estimated by the maximum likelihood estimation based on the peak-over-threthold (POT) developed by Pickands (1975), where the estimates depend on $x$ and the monotonicity is not ensured. The following result on the  convergence rate for $\gamma>0$ is a direct consequence of Theorem 8.1. in Smith (1987).

\begin{cor}
Suppose {\rm (i)} and the positive constants $c_1$ and $c_2$ exists such that $x=c_1 u$ and $\sqrt{N} u^{-\beta} \to c_2$ conditional on $N$ and $u$. Then, 
\begin{align*}
&\{\overline{F}(x)\}^{-2} \mathbb{E}[N({F}(x) - {\widehat{H}}(x))^2 | N,u] \\
\sim& B^2 \{ c_1^{-\beta} + c_2\alpha (\alpha+\beta)^{-1} (\alpha+\beta+1)^{-1} (\alpha-\beta, \alpha) \bm{\nu}\}^2 + 1+ \alpha^2 \bm{\nu}^{\mathsf{T}} \Sigma_0 \bm{\nu},
\end{align*}
where $\bm{\nu}^{\mathsf{T}} :=(c_1^{-1}-1, \alpha (\ln c_1+ c_1^{-1} -1) )$ and the asymptotic bias is zero if $\beta$ equals $0$ or $1$. For $x$ s.t. $u=o(x)$ and $N^{-1/2}\ln (1+u^{-1} x) \to 0$
\begin{align*}
&\{\overline{F}(x)\}^{-2} \mathbb{E}[N \{\ln (1+u^{-1} x)\}^{-2}({F}(x) - {\widehat{H}}(x))^2 | N,u] \\
\sim& B^2 \{c_2\alpha(\alpha+\beta)^{-1}(\alpha+\beta+1)^{-1} (\alpha+1)(1-\beta)\}^2 + \alpha^2(\alpha+1)^2.
\end{align*}
\end{cor}

The following result on the convergence rate for $\gamma=0$ is a direct consequence of Theorem Theorem 9.5. and 9.6. in Smith (1987).
\begin{cor}
Suppose {\rm (ii)} and a positive constant $c_3$ exists such that $x=u(1+c_3 (\kappa C)^{-1} u^{-\kappa})$ conditional on $N$ and $u$. If $(\kappa C)^{-1} (1-\kappa) \sqrt{N} u^{-\kappa} \to c_4$, then
\begin{align*}
\{\overline{F}(x)\}^{-2} \mathbb{E}[N({F}(x) - {\widehat{H}}(x))^2 | N,u] \sim 1+2c_3^2 -c_3^3 + 4^{-1} c_3^4.
\end{align*}
If $(\kappa C)^{-2} (1-\kappa)^2 \sqrt{N} u^{-2\kappa} \to c_5$, then
\begin{align*}
&\{\overline{F}(x)\}^{-2} \mathbb{E}[N({F}(x) - {\widehat{H}}(x))^2 | N,u] \\
\sim& c_3^2 c_5^2 \{-(4 -\kappa (1-\kappa)) + 2^{-1}c_3 (5+ 2\kappa (1-\kappa)) + 6^{-1} c_3^2 (2 +\kappa (1-\kappa)) -8^{-1} c_3^3 \}^2 \\
& + 1+2c_3^2 -c_3^3 + 4^{-1} c_3^4.
\end{align*}
\end{cor}

$u_0=F^{-1}(1-n^{-1}N_0)$ supposed in this study means $u_0=O((N_0/n)^{-\gamma})$ for the Hall class and $u_0= O(\{C^{-1} (-\ln (N_0/n))\}^{1/\kappa})$ for the Weibull class. Suppose $u_0=O(N_0^{1/(2\beta)})$ required in Corollary 3 {\rm (i)} and $u_0=O(N_0^{1/(2\kappa)})$ required in Corollary 3 {\rm (ii)}. Then, $N_0=O(n^{2\beta/(2\beta +\alpha)})$ and $u_0=O(n^{1/(2\beta+\alpha)})$ for the Hall class. $n=O(N \exp(C\sqrt{N}))$ lead by $u_0= O(\{C^{-1} (-\ln (N_0/n))\}^{1/\kappa})$ and $u_0=O(N_0^{1/(2\kappa)})$ for the Weibull class means there exists some slowly varying function $\ell^*$ such that $N\sim C^{-2} (\ln n)^2 \ell^*(n)$.

\subsection{The plug-in approach}

If we know that $F$ belongs to the Hall class in advance, we can estimate the tail probability by fitting without misspecification. The following plug-in type of the tail probability estimator for the Hall class is proposed by Hall and Weissman (1997)
$$F_{\widehat{\bm{\theta}}}(x):=1-\widehat{A}x^{-\widehat{\alpha}},$$
where the parameters $(A,\alpha)$ are estimated by the $r$ largest order statistics. The following result was obtained by Hall and Weissman (1997).

\begin{proposition}
Suppose {\rm (i)} and there exist a positive constant $c_6$ exists such that $\{(\ln x)/(\ln n)\} \to c_6$. Then, the optimal $r$ is given by $r_0\sim s n^{2\beta/(\alpha+2\beta)}$ for a positive constant $s$ and 
$$\mathbb{E}[(F(x) - F_{\widehat{\bm{\theta}}}(x))^2] \sim 2^{-1} s^{-1} A^2 \alpha^2 \beta^{-1} (\alpha+2\beta)^{-1} \{c_6(\alpha+2\beta) -1\}^2 n^{-2} (\ln n)^2,$$
where $s:= A^{2\beta/(\alpha+2\beta)} B^{-2\alpha/(\alpha+2\beta)}\{\alpha(\alpha+\beta)^2 /(2\beta^3)\}^{\alpha/(\alpha+2\beta)}$.
\end{proposition}

For the Weibull class of distributions estimating $\kappa$ called the Weibull-tail coefficient has been studied (e.g. Allouche et al. 2023). However, the parameter $C$ has been out of our target of estimation. We don't have any approach to assess $C$ to the best of our knowledge and to estimate the tail probability of Weibull class directly.

\section{Comparative study on tail probability estimators}

The asymptotic convergence rates of the estimators and numerical accuracy by a simulation study are investigated. The piecing-together estimator (PT) needs to determine $u$ and $N$ by some data-driven approach; however, we suppose $u$ and $N$ are nonrandom and $u\sim C_3 x$ and $N=O(n \overline{F}(u))$ following the study Smith (1987), where $0<C_3\le 1$ for the Hall class and $C_3=1$ for the Weibull class. Then, Corollary 3 requires $x\sim C_1 n^{\delta}$ for the Hall class and $x\sim C_2 (\ln n)^{1/\kappa}$ for the Weibull class, where 
$$\delta:=1/(2\beta+\alpha).$$
The convergence rate of the relative MSE of PT given by $\{\overline{F}(x)\}^{-2} (F(x) -\widehat{H}(x))^2$ is of order for the Hall class
$$n^{-1} x^{\alpha} = O(n^{\delta\alpha-1})$$
or for the Weibull class
$$n^{-1} \exp(C\{2x^{\kappa} -u^{\kappa}\}) = O(n^{-1} \exp( C x^{\kappa})) = O(n^{-1+C C_2^{\kappa}}).$$

The convergence rate of the nonparametric estimator (NE) with the optimal bandwidth is given in Corollary 1; however, 
$$n^{-1} \{f'(x)\}^{-2} f(x) =O(n^{-1} x^{\alpha+3}) \to 0$$ 
is required for the Hall class. That means $(\alpha+3)\delta<1 \Longleftrightarrow \beta>3/2$. The Weibull class always satisfies the condition of the bandwidth convergence. The optimal convergence rate of the relative MSE of NE is $n^{-1} x^{\alpha}$ for the Hall class or $n^{-1} \exp(Cx^{\kappa})$ for the Weibull class. That is the convergence rate is same as that of PT for these cases.

The convergence rate of the relative MSE of the plug-in estimator (PI) $F_{\widehat{\bm{\theta}}}(x)$ with the optimal $r=r_0$ for the Hall class is of order $n^{2\delta\alpha-2}(\ln n)^2$. The relative convergence rate of the different approaches are summarized in Table 1. The hyphens means the optimal bandwidth diverge, that is Corollary 1 is broken.
\begin{table}
\caption{The convergence rate of the relative MSE of the tail probability estimators at $x\sim C_1 n^{\delta}$ for the Hall class and $x\sim C_2 (\ln n)^{1/\kappa}$ for the Weibull class}
{\fontsize{10pt}{10pt}\selectfont
$$\begin{tabu}[c]{|c||ccc|}
  \hline
 & \text{NE} & \text{PT} & \text{PI } \\ \hline
\text{Hall class} & n^{\delta\alpha-1} & n^{\delta\alpha-1} & n^{2\delta\alpha-2} (\ln n)^2 \\
\text{Weibull class} & n^{-1+C C_2^{\kappa}} & n^{-1+C C_2^{\kappa}} & \text{--} \\ 

\hline
\end{tabu}$$
}
\end{table}

Comparing the convergence rates of the approaches in the Hall cases, PI outperforms the other approaches if $\delta\alpha<1$. However, it is always true since $\gamma\beta>0$. We see that the convergence rate becomes two times if we suppose the class of the underlying distribution the Hall class in advance and it is true.

If $C=C_2=1$ the relative convergence rate of NE and PT is $0$. If $C=1$, both approaches do not converge for $C_2 \ge 1$. Examples of the polynomial rates of the relative MSE of NE, PT and PI are provided in Table 2. For the Hall class the convergence rates depend on $\beta/\alpha$. As $\beta/\alpha$ gets large the convergence rates become faster. For the Weibull class the convergence rate becomes close to $-1$ as $\kappa$ gets large.

\begin{table}
\caption{Examples of the polynomial rates of the relative MSE of the tail probability estimators given in Table 1}
\begin{minipage}[t]{.5\textwidth}
{\fontsize{10pt}{10pt}\selectfont
$$\begin{tabu}[c]{|c|cc||ccc|}
  \hline
  \multicolumn{6}{|c|}{{\rm Burr}} \\ \hline
c,\ell & \alpha & \beta & {\rm NE} & {\rm PT} & {\rm PI} \\ \hline

 1/2,1/2 & 1/4 & 1/2 & \text{--} & -0.8 & -1.6 \\
 1,1/2 & 1/2 & 1 & \text{--} & -0.8 & -1.6 \\
 3,1/2 & 3/2 & 3 & -0.8 & -0.8 & -1.6 \\ \hline
 1/2,1 & 1/2 & 1/2 & \text{--} & -0.667 & -1.333 \\
 1,1 & 1 & 1 & \text{--} & -0.667 & -1.333 \\
 3,1 & 3 & 3 & -0.667 & -0.667 & -1.333 \\ \hline
 1/2,3 & 3/2 & 1/2 & \text{--} & -0.4 & -0.8 \\
 1,3 & 3 & 1 & \text{--} & -0.4 & -0.8 \\
 3,3 & 9 & 3 & -0.4 & -0.4 & -0.8 \\ \hline
 
\end{tabu}$$
}
 \end{minipage}
\begin{minipage}[t]{.5\textwidth}
{\fontsize{10pt}{10pt}\selectfont
$$\begin{tabu}[c]{|cc|c||cc|}
\hline
 \multicolumn{5}{|c|}{{\rm Weibull}} \\ \hline
\kappa& C & C_2 & {\rm NE} & {\rm PT} \\ \hline

1/2 & 1 & 1/5 & -0.553 & -0.553 \\
1 & 1 & 1/5 & -0.8 & -0.8 \\
3 & 1 & 1/5 & -0.992 & -0.992 \\
10 & 1 & 1/5 & -1.000 & -1.000 \\ \hline

1/2 & 1 & 1/3 & -0.423 & -0.423 \\
1 & 1 & 1/3 & -0.667 & -0.667 \\
3 & 1 & 1/3 & -0.963 & -0.963 \\
10 & 1 & 1/3 & -1.000 & -1.000 \\ \hline

1/2 & 1 & 1/2 & -0.293 & -0.293 \\
1 & 1 & 1/2 & -0.5 & -0.5 \\
3 & 1 & 1/2 & -0.875 & -0.875 \\
10 & 1 & 1/2 & -0.999 & -0.999 \\ \hline

\end{tabu}$$
}
\end{minipage}
\end{table}

\clearpage

\newgeometry{bottom=20pt}

\begin{landscape}
\begin{table}
\caption{Relative MSE values ($\times 100$) and sd values ($\times 100$) for the estimators with $n=2^{8}$ i.e. $x=C_1 2^{8\delta}$ or $n=2^{12}$ i.e. $x=C_1 2^{12\delta}$}
{\fontsize{6pt}{6pt}\selectfont
$$\begin{tabu}[c]{|c||cccccccc|cccccccc|cccccccc|}
   \hline
 {\rm Burr} & {\rm PI} & {\rm sd} & {\rm PT} & {\rm sd} & {\rm AL} & {\rm sd} & {\rm PB} & {\rm sd} & {\rm PI} & {\rm sd} & {\rm PT} & {\rm sd} & {\rm AL} & {\rm sd} & {\rm PB} & {\rm sd} & {\rm PI} & {\rm sd} & {\rm PT} & {\rm sd} & {\rm AL} & {\rm sd} & {\rm PB} & {\rm sd} \\  \hline\hline 
    
 c,\ell & \multicolumn{8}{c|}{n=2^{8}, C_1=1/2} & \multicolumn{8}{c|}{n=2^{8}, C_1=1} & \multicolumn{8}{c|}{n=2^{8}, C_1=2} \\ \hline
1/2,1/2 & 0.496 & 0.055 & 0.100 & 0.007 & 0.050 & 0.004 & 0.002 & 0.003 & 0.327 & 0.845 & 0.098 & 0.182 & 0.018 & 0.002 & 0.003 & 0.004 & 0.220 & 0.573 & 0.157 & 0.311 & 0.001 & 0.000 & 0.004 & 0.005 \\ 
  1,1/2 & 0.400 & 1.125 & 0.046 & 0.019 & 0.019 & 0.009 & 0.002 & 0.003 & 0.182 & 0.476 & 0.081 & 0.037 & 0.014 & 0.013 & 0.004 & 0.005 & 0.093 & 0.224 & 0.118 & 0.062 & 0.219 & 0.047 & 0.006 & 0.009 \\ 
  3,1/2 &  0.280 & 1.052 & 0.420 & 0.052 & 0.005 & 0.006 & 0.003 & 0.004 & 0.041 & 0.100 & 2.161 & 0.422 & 0.080 & 0.237 & 0.012 & 0.017 & 0.039 & 0.061 & 3.225 & 1.220 & 0.309 & 1.842 & 0.043 & 0.061 \\ 
  1/2,1 & 0.152 & 0.357 & 0.073 & 0.051 & 0.407 & 0.062 & 0.008 & 0.011 & 0.076 & 0.177 & 0.094 & 0.070 & 1.168 & 0.108 & 0.012 & 0.016 & 0.043 & 0.095 & 0.118 & 0.096 & 2.834 & 0.337 & 0.016 & 0.021 \\ 
  1,1 & 0.078 & 0.290 & 0.392 & 0.153 & 0.463 & 0.465 & 0.010 & 0.016 &  0.031 & 0.072 & 0.622 & 0.306 & 1.354 & 2.185 & 0.022 & 0.033 & 0.037 & 0.054 & 0.832 & 0.523 & 3.455 & 8.035 & 0.046 & 0.068 \\ 
  3,1 & 0.058 & 0.156 & 4.138 & 1.365 & 0.013 & 0.018 & 0.012 & 0.018 & 0.076 & 0.101 & 26.62 & 12.29 & 0.101 & 0.159 & 0.107 & 0.169 & 0.457 & 0.847 & 27.31 & 38.40 & 0.800 & 1.270 & 2.458 & 3.173 \\ 
  1/2,3 & 0.099 & 0.155 & 1.298 & 1.511 & 2.791 & 21.61 & 0.398 & 0.446 & 0.297 & 0.486 & 1.829 & 2.783 & 4.212 & 47.79 & 6.314 & 2.568 & 0.957 & 1.513 & 2.797 & 6.108 & 4.208 & 64.13 & 136.8 & 61.48 \\ 
  1,3 & 0.070 & 0.097 & 5.560 & 4.428 & 0.128 & 0.182 & 0.141 & 0.243 & 0.607 & 0.952 & 9.385 & 15.33 & 0.618 & 0.974 & 2.367 & 2.032 & 5.748 & 8.872 & 11.81 & 54.39 & 3.549 & 7.972 & 1546 & 1428 \\ 
  3,3 & 0.094 & 0.232 & 0.012 & 0.016 & 0.014 & 0.020 & 0.013 & 0.019 & 1.234 & 2.009 & 1039 & 1641 & 1.047 & 1.846 & 2.070 & 2.966 & 676.1 & 1662 & 2 \times 10^4 & 2 \times 10^5 & 166.0 & 4458 & 6 \times 10^4 & 8 \times 10^4 \\ 
   
 \hline\hline

 c,\ell & \multicolumn{8}{c|}{n=2^{12}, C_1=1/2} & \multicolumn{8}{c|}{n=2^{12}, C_1=1} & \multicolumn{8}{c|}{n=2^{12}, C_1=2} \\ \hline
1/2,1/2 & 0.125 & 0.175 & 0.302 & 0.231 & 0.023 & 0.000 & 0.000 & 0.000 & 0.081 & 0.123 & 0.560 & 0.386 & 0.002 & 0.000 & 0.000 & 0.000 & 0.055 & 0.088 & 0.918 & 0.657 & 0.009 & 0.000 & 0.000 & 0.000 \\ 
  1,1/2 & 0.062 & 0.098 & 0.073 & 0.009 & 0.001 & 0.001 & 0.000 & 0.000 & 0.030 & 0.053 & 0.107 & 0.014 & 0.140 & 0.014 & 0.000 & 0.001 & 0.017 & 0.030 & 0.134 & 0.020 & 0.763 & 0.038 & 0.001 & 0.001 \\ 
  3,1/2 & 0.018 & 0.034 & 1.616 & 0.070 & 0.006 & 0.042 & 0.001 & 0.001 & 0.005 & 0.008 & 2.999 & 0.236 & 0.034 & 0.507 & 0.002 & 0.003 & 0.005 & 0.006 & 3.318 & 0.473 & 0.141 & 3.730 & 0.006 & 0.008 \\ 
  1/2,1 & 0.044 & 0.061 & 0.095 & 0.125 & 1.945 & 0.069 & 0.001 & 0.001 & 0.020 & 0.033 & 0.108 & 0.033 & 4.716 & 0.109 & 0.001 & 0.002 & 0.010 & 0.017 & 0.121 & 0.033 & 10.64 & 0.217 & 0.002 & 0.003 \\ 
  1,1 & 0.009 & 0.014 & 0.647 & 0.088 & 1.148 & 2.751 & 0.002 & 0.003 & 0.004 & 0.006 & 0.806 & 0.145 & 2.665 & 9.474 & 0.004 & 0.006 & 0.006 & 0.009 & 0.917 & 0.234 & 6.398 & 31.21 & 0.008 & 0.011 \\ 
  3,1 & 0.004 & 0.006 & 23.87 & 1.353 & 0.005 & 0.007 & 0.021 & 0.031 & 0.004 & 0.006 & 44.41 & 6.977 & 0.037 & 0.053 & 0.035 & 0.053 & 0.084 & 0.139 & 42.40 & 28.44 & 0.299 & 0.480 & 0.336 & 0.562 \\ 
  1/2,3 & 0.040 & 0.059 & 1.813 & 0.891 & 1.445 & 44.19 & 0.049 & 0.075 & 0.160 & 0.196 & 2.161 & 1.647 & 0.457 & 10.89 & 0.175 & 0.243 & 0.517 & 0.553 & 2.524 & 3.136 & 0.268 & 0.362 & 19.84 & 12.11 \\  
  1,3 & 0.040 & 0.054 & 15.70 & 4.180 & 0.044 & 0.063 &  0.047 & 0.069 & 0.508 & 0.526 & 20.85 & 18.02 & 0.286 & 0.377 & 0.364 & 0.587 & 3.441 & 3.271 & 18.06 & 48.22 & 1.963 & 3.512 & 1226 & 1383 \\ 
  3,3 & 0.005 & 0.007 & 0.029 & 0.018 & 0.004 & 0.007 & 0.004 & 0.006 & 1.209 & 1.186 & 2 \times 10^4 & 10^4 & 0.625 & 0.921 & 0.843 & 1.413 & 145.5 & 182.1 & 3 \times 10^4 & 2 \times 10^5 & 350.9 & 5524 & 2 \times 10^6 & 2 \times 10^6 \\ 

\hline

\end{tabu}$$
}
\end{table}

\begin{table}[h]
\caption{Relative MSE values ($\times 100$) and sd values ($\times 100$) for the estimators with $n=2^{8}$  i.e. $x=C_2 (\ln 8)^{1/\kappa}$ or $n=2^{12}$ i.e. $x=C_2 (\ln 12)^{1/\kappa}$}
{\fontsize{6pt}{6pt}\selectfont
$$\begin{tabu}[c]{|c||cccccc|cccccc|cccccc|}
   \hline
{\rm Weibull} & {\rm PT} & {\rm sd} & {\rm AL} & {\rm sd} & {\rm PB} & {\rm sd} & {\rm PT} & {\rm sd} & {\rm AL} & {\rm sd} & {\rm PB} & {\rm sd} & {\rm PT} & {\rm sd} & {\rm AL} & {\rm sd} & {\rm PB} & {\rm sd} \\  \hline\hline 
    
\kappa, C & \multicolumn{6}{c|}{n=2^{8}, C_2=1/5} & \multicolumn{6}{c|}{n=2^{8}, C_2=1/3} & \multicolumn{6}{c|}{n=2^{8}, C_2=1/2}\\ \hline
1/2,1 & 0.164 & 3.784 & 0.041 & 0.062 & 2.747 & 0.472 & 0.096 & 0.128 & 0.088 & 0.138 & 69.41 & 1.532 & 0.200 & 0.279 & 0.189 & 0.304 & 461.7 & 2.474 \\ 
  1,1 & 0.008 & 0.010 & 0.008 & 0.010 & 0.008 & 0.012 & 0.020 & 0.026 & 0.019 & 0.027 & 0.018 & 0.026 & 0.064 & 0.095 & 0.052 & 0.073 & 0.054 & 0.086 \\ 
  3,1 & 0.000 & 0.000 & 0.000 & 0.000 & 0.000 & 0.000 & 0.001 & 0.001 & 0.001 & 0.001 & 0.001 & 0.001 & 0.004 & 0.005 & 0.003 & 0.004 & 0.004 & 0.006 \\ 
  10,1 & 0.000 & 0.000 & 0.000 & 0.000 & 0.000 & 0.000 & 0.000 & 0.000 & 0.000 & 0.000 & 0.000 & 0.000 & 0.000 & 0.000 & 0.000 & 0.000 & 0.000 & 0.000 \\ 
  
  \hline\hline
  
  \kappa, C & \multicolumn{6}{c|}{n=2^{12}, C_2=1/5} & \multicolumn{6}{c|}{n=2^{12}, C_2=1/3} & \multicolumn{6}{c|}{n=2^{12}, C_2=1/2}\\ \hline
1/2,1 & 0.011 & 0.015 & 0.010 & 0.014 & 240.5 & 0.699 & 0.031 & 0.046 & 0.030 & 0.045 & 3038 & 1.681 & 0.097 & 0.146 & 0.093 & 0.128 & 29460 & 4.633 \\ 
  1,1 & 0.008 & 0.072 & 0.001 & 0.001 & 0.001 & 0.001 & 0.004 & 0.023 & 0.004 & 0.005 & 0.009 & 0.009  & 0.016 & 0.023 & 0.015 & 0.021 & 1.784 & 0.281 \\ 
  3,1 & 0.000 & 0.000 & 0.000 & 0.000 & 0.000 & 0.000 & 0.000 & 0.000 & 0.000 & 0.000 & 0.000 & 0.000 & 0.000 & 0.001 & 0.000 & 0.001 & 0.000 & 0.001 \\ 
  10,1 & 0.000 & 0.000 & 0.000 & 0.000 & 0.000 & 0.000 & 0.000 & 0.000 & 0.000 & 0.000 & 0.000 & 0.000 & 0.000 & 0.000 & 0.000 & 0.000 & 0.000 & 0.000 \\ 
  
     \hline
\end{tabu}$$
}
\end{table}
\end{landscape}

\restoregeometry

Simulation studies were conducted on the relative MSE of randomly generated 1000 datasets following each Burr distribution. In this study $(A, \alpha)$ are estimated by the $r=n^{2/3}$ largest order statistics with the supposition $\alpha=\beta=1$. Table 3 reports that of PI 
$$\{\overline{F}(x)\}^{-2} \{F(x) - F_{\widehat{\bm{\theta}}}(x)\}^2,$$
PT and NE with the Altman and L\'{e}ger (1995)'s bandwidth estimator (AL) provided by \texttt{aLbw} function in the \texttt{kerdiest} package in R. PB corresponds that of NE with the pointwise optimal bandwidth estimator provided in Section 2. The kernel functions of NE are the Gaussian. 

The Burr distribution belongs to the Hall class with $\alpha=\beta\ell$. That means the performance of the estimators theoretically depends on $\ell$. To be exact, as $\ell$ gets large the convergence rates become slower. $u=0.5x$ (i.e. $C_3=0.5$) for the Burr distribution and $u=0.99x$ (i.e. $C_3=0.99$) for the Weibull distribution are chosen for PT in this simulation study. 

The table shows the rough coincidence with the theoretical property; however, the performance also depends on $C_1$. When at least one of $\beta, \ell$ or $C_1$ is small, PB is sometimes best. However, PB is numerically unstable, especially for $\beta=3, \ell=3, C_1=2$. Additional simulation studies demonstrate that this mainly comes from the bandwidth $h_{A_0, \alpha_0}$ with the exact values of the underlying distribution parameters $(A_0, \alpha_0)$ is considerably larger than that of Altman and L\'{e}ger (1995)'s bandwidth, which causes the under-smoothing. When $\beta=3, \ell=3, C_1=2$, $\overline{F}(x)$ is about $10^{-5}$ for $n=2^8$, $10^{-6}$ for $n=2^{12}$ and $10^{-10}$ for $n=2^{24}$. Nevertheless the estimated values by PB are about $0.03$, $0.01$ and $0.002$, respectively. When $\beta=3, \ell=3, C_1=2$, PI tends to work best. 

The next study is on the relative MSE of PT, AL and PB for the Weibull class, which is shown in Table 4. 1000 datasets were generated for each of the Weibull distribution.  Simulation studies for PI (the plug-in estimator for the Hall class) demonstrated that PI cannot be applied for the Weibull class at all. Table 4 shows that the estimators become precise as $\kappa$ gets large. On the whole AL tends to slightly outperform PT. The bandwidth of PB is also given by $h_{\widehat{A},\widehat{\alpha}}$, and so the case is considered to be misspecified. PB does not work especially for $\kappa=1/2$. 

Comparing the MSE values, we see that the estimators are not one-sided. In the sense of the numerical stability AL tends to be best among them in the setting we cannot suppose the underlying distribution belongs to the Hall class.

\section{Mean excess function estimation}
\subsection{The kernel type estimator as a nonparametric approach}

This section considers estimation of MEF. The kernel type estimator is given by
$$\widehat{e}(u) := \int_{0}^{\infty} x \frac{\widehat{f}(x+u)}{1-\widehat{F}(u)} {\rm d}x,$$
where $\widehat{f}(x)$ is the derivative of the kernel distribution estimator $\widehat{F}(x)$. 

Mugdadi and Teweldemedhin (2013) proposes two different type of MEF estimators; however, this study employs the above kernel type MEF estimator referred to as KE. Asymptotic properties of the kernel type estimator are provided in Guillamon et al. (1998) (see also for censored data Chaubey and Sen 2008 and for length-biased data Zamini et al. 2023). 

\begin{theorem}
Suppose $F$ belongs to the Hall class with $\alpha>1$ or the Weibull class. Then, if $\int z^2 w(z) {\rm d}z <\infty$ and $\vert\int z W(z) w(z) {\rm d}z\vert < \infty$
\begin{align*}
\mathbb{E}[(\widehat{e}(u) - e(u))^2] \sim & h^4 \left( \int z^2 w\left(z\right) {\rm d}z \right)^2
\begin{dcases}
\alpha^2 (\alpha-1)^{-2} u^{-2} \\
o(u^{4\kappa -2})
\end{dcases}
\\
+ &\frac{1}{n}
\begin{dcases}
2A^{-1}(\alpha -1)^{-2} u^{\alpha+2} \{(\alpha-2)^{-1} (\alpha-1) - \alpha h u^{-1} \int z W(z) w(z) {\rm d}z\}\\
o(u^{2} \exp(Cu^{\kappa})).
\end{dcases}
\end{align*}
If the class is the Hall, the asymptotically optimal bandwidth is
$$h_{A,\alpha}^*=n^{-1/3} u^{1+(\alpha/3)} \left\{2^{-1} A^{-1} \alpha^{-1} \frac{\int z W(z) w(z) {\rm d}z}{(\int z^2 w\left(z\right) {\rm d}z)^2}\right\}^{1/3}$$
\end{theorem}

The requirement $h_{A,\alpha}^* \to 0$ is $n^{-1}u^{\alpha+3} \to 0$ meaning $u=o(n^{1/(\alpha+3)})$. The convergence rate of the MSE with the bandwidth $h_{A,\alpha}^*$ is $n^{-1}u^{\alpha+2}$, which converges to zero if $h_{A,\alpha}^* \to 0$. 

\begin{cor}
Suppose $F$ belongs to {\rm (i)} the Hall class, $r=O(n^{2\beta/(2\beta+\alpha)})$ and $u=o(\exp(n^{\beta/(2\beta+\alpha)}))$. Then, under the assumptions of Theorem 1, 
$$(h_{A,\alpha}^*)^{-1} (h_{\widehat{A},\widehat{\alpha}}^* -h_{A,\alpha}^*) = O_P(n^{-\beta/(2\beta+\alpha)} (\ln x + \ln n)),$$
where $\widehat{A}$ and $\widehat{\alpha}$ are estimator consists of the $r$ largest order statistics.
\end{cor}

\subsection{Fitting estimator to GPD}

The Pickands-Balkema-De Haan theorem in the extreme value theory is applicable to the mean excess function $F_u(x)$. The generalized pareto distribution estimator $H_{\widehat{\bm{\gamma}}}$ yields the following MEF estimator
$$e_{\widehat{\bm{\gamma}}}(u) := \frac{\widehat{c}_u}{1-\widehat{\gamma}}$$
which is considered as a parametrically fitting estimator (PE). Let $N$ be the number of $X_1,X_2,\cdots X_n$ exceeding $u$ and ${N^*}:=n(1 -F(u))$. Applying the asymptotic normality of the scaled parameter vector
$$\widehat{\bm{\gamma}}^*=\begin{pmatrix}
1 & \\
 & c_u^{-1} 
\end{pmatrix}
\widehat{\bm{\gamma}}$$ 
(Smith 1987), we see that the following corollary is a consequence of the delta method.

\begin{cor}
Suppose that $F$ belongs to the Hall class with $\alpha>1$ or the Weibull class and that $^{\exists} \lambda \in \mathbb{R}$ s.t. $\lambda_n \to \lambda$, where
$$\lambda_n := \sqrt{n} \times 
\begin{dcases}
-A^{1/2} B \beta (\alpha+\beta)^{-1} u^{-(\alpha/2)-\beta} ~~~ &{\rm for ~ the ~ Hall ~ class} \\
(C \kappa)^{-2} u^{-2\kappa} \exp(-2^{-1} Cu^{\kappa})~~~ &{\rm for ~ the ~ Weibull ~ class}.
\end{dcases}
$$
Then, $\sqrt{N^*} (e(u) - e_{\widehat{\bm{\gamma}}}(u))$ converges in the normal distribution with the asymptotic mean 
$$\begin{dcases}
\lambda_n e(u) (1+\gamma)\{1 + \gamma (1-\gamma)^{-1} (1-\beta) (\alpha+\beta+1)^{-1}\}~~~ &{\rm for ~ the ~ Hall ~ class} \\
\lambda_n u (\kappa-1) ~~~ &{\rm for ~ the ~ Weibull ~ class}.
\end{dcases}$$
and the asymptotic variance
$$\begin{dcases}
e^2(u) (1+\gamma) (1-\gamma)^2 (2\gamma^2 - \gamma +1) ~~~ &{\rm for ~ the ~ Hall ~ class} \\
C^{-2} \kappa^{-2} u^{2(1-\kappa)} ~~~ &{\rm for ~ the ~ Weibull ~ class},
\end{dcases}$$
hence
\begin{align*}
&\mathbb{E}[(e(u) - e_{\widehat{\bm{\gamma}}}(u))^2] \\
\sim& (N^*)^{-1}
\begin{dcases}
\lambda_n^2 e^2(u) (1+\gamma)^2\{1 + \gamma (1-\gamma)^{-1} (1-\beta) (\alpha+\beta+1)^{-1}\}^2 \\
~~~~~~ + e^2(u) (1+\gamma) (1-\gamma)^2 (2\gamma^2 - \gamma +1) ~~~ &{\rm for ~ the ~ Hall ~ class} \\
\lambda_n^2 u^2 (\kappa-1)^2 + C^{-2} \kappa^{-2} u^{2(1-\kappa)} ~~~ &{\rm for ~ the ~ Weibull ~ class}.
\end{dcases}
\end{align*}
\end{cor}

The convergence rate of MSE for the Hall class $(N^*)^{-1} e^2(u) = O(n^{-1} u^{\alpha+2})$ is same as that of KE. However, the condition is different. KE requires $n^{-1} u^{\alpha+3} \to 0$, while $n u^{-\alpha-2\beta} \to 0$

\subsection{The plug-in approach}

If we know that $F$ belongs to the Hall class in advance, the following plug-in type of the estimator (PT) can be applied
$$e_{\widehat{\bm{\theta}}}(u):= \frac{u}{\widehat{\alpha} -1} = \int_{0}^{\infty} x \widehat{\alpha} u^{\widehat{\alpha}} (x+u)^{-\widehat{\alpha}-1} {\rm d}x,$$
where $\alpha$ is estimated by the $r$ largest order statistics. The following result is a consequence of Hall (1982).

\begin{cor}
If $F$ belongs to {\rm (i)} the Hall class and $r\sim sn^{2\beta/(2\beta+\alpha)}$ for a positive constant $s$, 
$u^{-1} n^{(1/2)-(\alpha/(2\beta+\alpha))/2} (e(u) - e_{\widehat{\bm{\theta}}}(u))$ converges in the normal distribution with the asymptotic mean 
$$-\frac{\alpha}{(\alpha-1)^2} A^{-\beta/\alpha}B\beta (\alpha+\beta)^{-1} s^{\beta/\alpha}$$
and the asymptotic variance $\alpha^2 (\alpha-1)^{-2} s^{-1}$, hence,
$$\mathbb{E}[(e_{\widehat{\bm{\theta}}}(u) -e(u))^2] \sim u^2 n^{-1+(\alpha/(2\beta+\alpha))} \frac{\alpha^2}{(\alpha-1)^4} (A^{-2\beta/\alpha}B^2\beta^2 (\alpha+\beta)^{-2} s^{2\beta/\alpha} + s^{-1}).$$
\end{cor}

Interestingly, different from that of the tail probability estimator, the convergence rate $u^2 n^{-1+(\alpha/(2\beta+\alpha))}$ is not necessarily faster than that of KE and PE, which is true if $n =o(u^{\alpha+2\beta})$. 

\section{Comparative study on mean excess function estimators}

Table 5 shows the convergence rate of $\mathbb{E}[(\widehat{e}(u) - e(u))^2]$ with $h=h_{A,\alpha}^*$ for the Hall class of distributions, that of PE and that of PT. The hyphen in Table 5 means that the assumption of the obtained convergence rate is broken. Table 5 shows that the convergence rate of NE is fastest in the settings as long as $h_{A,\alpha}^* \to 0$. For large $u$ only PI has the consistency.

\begin{table}
\caption{Examples of the polynomial rates of the MSE of the mean excess function estimators}
{\fontsize{8pt}{8pt}\selectfont
$$\begin{tabu}[c]{|ccc||ccc|ccc|ccc|ccc|}
  \hline
  \multicolumn{3}{|c|}{{\rm Hall}} & \multicolumn{3}{c|}{u=n^{1/16}} & \multicolumn{3}{c|}{u=n^{1/8}} & \multicolumn{3}{c|}{u=n^{1/4}} & \multicolumn{3}{c|}{u=n^{3/8}} \\ \hline
c,\ell & \alpha & \beta & {\rm NE} & {\rm PT} & {\rm PI} & {\rm NE} & {\rm PT} & {\rm PI} & {\rm NE} & {\rm PT} & {\rm PI} & {\rm NE} & {\rm PT} & {\rm PI} \\ \hline

 1/2,1/2 & 1/4 & 1/2 & -0.859 & \text{--} & -0.675 & -0.719 & \text{--} & -0.55 & -0.438 & \text{--} & -0.3 & \text{--} & \text{--} & -0.05 \\
 1,1/2 & 1/2 & 1 & -0.844 & \text{--} & -0.675 & -0.688 & \text{--} & -0.55 & -0.375 & \text{--} & -0.3 & \text{--} & \text{--} & -0.05 \\
 3,1/2 & 3/2 & 3 & -0.781 & \text{--} & -0.675 & -0.563 & \text{--} & -0.55 & \text{--} & -0.125 & -0.3 & \text{--} & 0.313 & -0.05 \\
 1/2,1 & 1/2 & 1/2 & -0.844 & \text{--} & -0.542 & -0.688 & \text{--} & -0.417 & -0.375 & \text{--} & -0.167 & \text{--} & \text{--} & 0.083 \\
 1,1 & 1 & 1 & -0.813 & \text{--} & -0.542 & -0.625 & \text{--} & -0.417 & -0.25 & \text{--} & -0.167 & \text{--} & 0.125 & 0.083 \\
 3,1 & 3 & 3 & -0.688 & \text{--} & -0.542 & -0.375 & -0.375 & -0.417 & \text{--} & 0.25 & -0.167 & \text{--} & 0.875 & 0.083 \\
 1/2,3 & 3/2 & 1/2 & -0.781 & \text{--} & -0.275 & -0.563 & \text{--} & -0.15 & \text{--} & \text{--} & 0.1 & \text{--} & \text{--} & 0.35 \\
 1,3 & 3 & 1 & -0.688 & \text{--} & -0.275 & -0.375 & \text{--} & -0.15 & \text{--} & 0.25 & 0.1 & \text{--} & 0.875 & 0.35 \\
 3,3 & 9 & 3 & -0.313 & \text{--} & -0.275 & \text{--} & 0.375 & -0.15 & \text{--} & 1.75 & 0.1 & \text{--} & 3.125 & 0.35 \\
 
 \hline
 
\end{tabu}$$
}
\end{table}

This section surveys the following the relative MSE in finite-sample cases of PE
$$\{e(u)\}^{-2} \{e_{\widetilde{\bm{\gamma}}}(u) -e(u)\}^2,$$
that of PI $e_{\widehat{\bm{\theta}}}(u)$ and NE $\widehat{e}(u)$, where
$$e_{\widetilde{\bm{\gamma}}}(u) := \frac{\widetilde{c}_u}{1-\widetilde{\gamma}}$$
and $\widetilde{\bm{\gamma}}:=(\widetilde{\gamma}, \widetilde{c}_u)$ is the maximum likelihood estimator under the constraint $\widetilde{\gamma}<1$. $(A, \alpha)$ are estimated by the $r=n^{2/3}$ largest order statistics. The kernel of NE $\widehat{e}(u)$ is set to the Gaussian, then it holds that
$$\{1-\widehat{F}(u)\}^{-1} n^{-1} \sum_{i=1}^n \left\{(X_i - u) \Phi\left(\frac{X_i -u}{h}\right) + h \phi\left(\frac{X_i -u}{h}\right)\right\}.$$ 
The bandwidth is the Altman and L\'{e}ger (1995)'s or $h_{\widehat{A},\widehat{\alpha}}^*$ provided in Corollary 2. 

We simulated the MISE values 10000 times, where Tables 6--7 show the mean values and their standard deviation (sd), where the sample sizes were $(n=)2^{8}$ {or} $2^{12}$. $u$ is 0.9, 0.95, 0.99, 0.995 \% quantiles of the underlying distributions, which is supposed to be the Burr or the Weibull.

\newgeometry{bottom=20pt}

\begin{landscape}
\begin{table}
\caption{Relative MSE values ($\times 100$) and sd values ($\times 100$) for the estimators with $n=2^{8}$ or $n=2^{12}$}
{\fontsize{5pt}{5pt}\selectfont
$$\begin{tabu}[c]{|c||cccccccc|cccccccc|cccccccc|}
   \hline
 {\rm Burr} & {\rm PI} & {\rm sd} & {\rm PT} & {\rm sd} & {\rm AL} & {\rm sd} & {\rm PB} & {\rm sd} & {\rm PI} & {\rm sd} & {\rm PT} & {\rm sd} & {\rm AL} & {\rm sd} & {\rm PB} & {\rm sd} & {\rm PI} & {\rm sd} & {\rm PT} & {\rm sd} & {\rm AL} & {\rm sd} & {\rm PB} & {\rm sd} \\  \hline\hline 
    
 c,\ell & \multicolumn{8}{c|}{n=2^{8}, C_1=1/2} & \multicolumn{8}{c|}{n=2^{8}, C_1=1} & \multicolumn{8}{c|}{n=2^{8}, C_1=2} \\ \hline
 
 3,1/2 & 7.448 & 141.9 & 2\times10^{19} & 8\times10^{20} & 0.277 & 2.261 & 1.151 & 12.67 & 7.448 & 141.9 & 5\times10^{16} & 10^{18} & 0.281 & 0.786 & 1.696 & 18.84 & 7.448 & 141.9 & 10^{18} & 3 \times10^{19} & \infty & \text{--} & 2.052 & 22.18 \\ 
  3,1 & 0.116 & 0.238 & 5\times10^{13} & 2\times10^{15} & 0.099 & 0.222 & 0.084 & 0.185 & 0.116 & 0.238 & 3\times10^{15} & 7\times10^{16} & 0.180 & 0.485 & 0.135 & 0.304 & 0.116 & 0.238 & 3\times10^{17} & 7\times10^{18} & \infty & \text{--} & 0.320 & 0.658 \\ 
  1/2,3 & 4924 & 76810 & 5\times10^{20} & 2\times10^{22} & 2.222 & 37.74 & 8.091 & 73.16 & 4924 & 76810 & 9\times10^{24} & 3\times10^{26} & 0.863 & 7.271 & 5.428 & 58.64 & 4924 & 76810 & 10^{19} & 4\times10^{20} & \infty & \text{--} & 2.798 & 46.10 \\ 
  1,3 & 5.213 & 9.553 & 8\times10^{14} & 2\times10^{16} & 1.046 & 1.556 & 0.975 & 1.554 & 5.213 & 9.553 & 5\times10^{14} & 6\times10^{15} & 0.784 & 1.904 & 0.610 & 1.633 & 5.213 & 9.553 & 9\times10^{17} & 3\times10^{19} & \infty & \text{--} & 0.445 & 1.689 \\ 
  3,3 & 0.914 & 0.659 & 0.565 & 0.621 & 0.499 & 0.473 & 0.497 & 0.417 & 0.914 & 0.659 & 4\times10^{17} & 10^{19} & 0.348 & 0.542 & 0.286 & 0.377 & 0.914 & 0.659 & 5\times10^{18} & 8\times10^{19} & 0.691 & 2.921 & 0.106 & 0.393 \\
   
 \hline\hline

 c,\ell & \multicolumn{8}{c|}{n=2^{12}, C_1=1/2} & \multicolumn{8}{c|}{n=2^{12}, C_1=1} & \multicolumn{8}{c|}{n=2^{12}, C_1=2} \\ \hline
 
 3,1/2 & 0.039 & 0.065 & 0.102 & 1.026 & 1.588 & 48.42 & 5.377 & 164.5 & 0.039 & 0.065 & 5\times10^9 & 10^{11} & 0.673 & 18.540 & 8.596 & 263.2 & 0.039 & 0.065 & 2\times10^{16} & 5\times10^{17} & 0.371 & 2.482 & 27.98 & 858.4 \\ 
  3,1 & 0.012 & 0.018 & 0.011 & 0.017 & 0.011 & 0.019 & 0.011 & 0.018 & 0.012 & 0.018 & 0.016 & 0.028 & 0.015 & 0.029 & 0.014 & 0.027 & 0.012 & 0.018 & 4\times10^6 & 10^8 & 0.062 & 0.128 & 0.055 & 0.105 \\ 
  1/2,3 & 3\times10^6 & 8\times10^7 & 5\times10^{13} & 2\times10^{15} & 2.325 & 2.089 & 4.455 & 18.13 & 3\times10^6 & 8\times10^7 & 3\times10^{12} & 5\times10^{13} & 1.109 & 1.612 & 2.498 & 14.31 & 3\times10^6 & 8\times10^7 & 3\times10^{19} & 10^{21} & 0.849 & 3.399 & 1.922 & 16.17 \\ 
  1,3 & 0.756 & 0.369 & 0.793 & 0.278 & 0.789 & 0.338 & 0.781 & 0.335 & 0.756 & 0.369 & 0.391 & 0.252 & 0.388 & 0.340 & 0.375 & 0.330 & 0.756 & 0.369 & 2\times10^{12} & 7\times10^{13} & 0.226 & 1.094 & 0.191 & 1.035 \\ 
  3,3 & 0.282 & 0.110 & 0.425 & 0.108 & 0.425 & 0.104 & 0.426 & 0.102 & 0.282 & 0.110 & 0.211 & 0.096 & 0.210 & 0.092 & 0.210 & 0.088 & 0.282 & 0.110 & 0.440 & 10.77 & 0.090 & 0.123 & 0.065 & 0.084 \\ 

\hline

\end{tabu}$$
}
\end{table}

\begin{table}[h]
\caption{Relative MSE values ($\times 100$) and sd values ($\times 100$) for the estimators with $n=2^{8}$  i.e. $x=C_2 (\ln 8)^{1/\kappa}$ or $n=2^{12}$ i.e. $x=C_2 (\ln 12)^{1/\kappa}$}
{\fontsize{5pt}{5pt}\selectfont
$$\begin{tabu}[c]{|c||cccccccc|cccccccc|cccccccc|}
   \hline
{\rm Weibull} & {\rm PI} & {\rm sd} & {\rm PT} & {\rm sd} & {\rm AL} & {\rm sd} & {\rm PB} & {\rm sd} & {\rm PI} & {\rm sd} & {\rm PT} & {\rm sd} & {\rm AL} & {\rm sd} & {\rm PB} & {\rm sd} & {\rm PI} & {\rm sd} & {\rm PT} & {\rm sd} & {\rm AL} & {\rm sd} & {\rm PB} & {\rm sd} \\  \hline\hline 
    
\kappa, C & \multicolumn{8}{c|}{n=2^{8}, C_2=1/5} & \multicolumn{8}{c|}{n=2^{8}, C_2=1/3} & \multicolumn{8}{c|}{n=2^{8}, C_2=1/2}\\ \hline

1/2,1 & 636.6 & 10400 & 6\times10^9 & 2\times10^{11} & 0.171 & 0.114 & 0.176 & 0.113 & 636.6 & 10400 & 10^{15} & 3\times10^{16} & 0.336 & 0.163 & 0.352 & 0.152 & 636.6 & 10400 & 2\times10^{16} & 4\times10^{17} & \infty & \text{--} & 0.650 & 0.164 \\ 
  1,1 & 0.158 & 0.252 & 0.047 & 0.069 & 0.045 & 0.055 & 0.041 & 0.049 & 0.158 & 0.252 & 3\times10^{11} & 8\times10^{12} & 0.150 & 0.119 & 0.140 & 0.099 & 0.158 & 0.252 & 4\times10^{18} & 10^{20} & \infty & \text{--} & 0.354 & 0.126 \\ 
  3,1 & 0.057 & 0.079 & 0.031 & 0.043 & 0.018 & 0.025 & 0.014 & 0.020 & 0.057 & 0.079 & 10^{10} & 4\times10^{11} & 0.072 & 0.067 & 0.021 & 0.021 & 0.057 & 0.079 & 7\times10^{15} & 10^{17} & 0.297 & 0.181 & 0.031 & 0.016 \\ 
  10,1 & 0.050 & 0.072 & 2\times10^{9} & 4\times10^{10} & 0.018 & 0.028 & 0.243 & 0.091 & 0.050 & 0.072 & 0.095 & 0.397 & 0.036 & 0.039 & 0.313 & 0.074 & 0.050 & 0.072 & 10^{17} & 2\times10^{18} & 0.188 & 0.117 & 0.872 & 0.244 \\ 
  
  \hline\hline
  
  \kappa, C & \multicolumn{8}{c|}{n=2^{12}, C_2=1/5} & \multicolumn{8}{c|}{n=2^{12}, C_2=1/3} & \multicolumn{8}{c|}{n=2^{12}, C_2=1/2}\\ \hline

1/2,1 & 0.142 & 0.055 & 0.138 & 0.030 & 0.142 & 0.029 & 0.143 & 0.029 & 0.142 & 0.055 & 0.304 & 0.044 & 0.308 & 0.043 & 0.310 & 0.042 & 0.142 & 0.055 & 0.540 & 0.071 & 0.545 & 0.069 & 0.552 & 0.064 \\ 
  1,1 & 0.056 & 0.026 & 0.019 & 0.011 & 0.019 & 0.011 & 0.019 & 0.011 & 0.056 & 0.026 & 0.113 & 0.031 & 0.113 & 0.030 & 0.113 & 0.029 & 0.056 & 0.026 & 0.326 & 0.076 & 0.328 & 0.074 & 0.327 & 0.061 \\ 
  3,1 & 0.036 & 0.017 & 0.005 & 0.005 & 0.004 & 0.004 & 0.003 & 0.004 & 0.036 & 0.017 & 0.062 & 0.023 & 0.060 & 0.020 & 0.053 & 0.018 & 0.036 & 0.017 & 0.239 & 0.069 & 0.238 & 0.061 & 0.185 & 0.036 \\ 
  10,1 & 0.032 & 0.016 & 0.003 & 0.004 & 0.002 & 0.002 & 0.002 & 0.002 & 0.032 & 0.016 & 0.051 & 0.020 & 0.044 & 0.017 & 0.010 & 0.006 & 0.032 & 0.016 & 0.268 & 1.779 & 0.202 & 0.057 & 0.016 & 0.004 \\ 

     \hline
\end{tabu}$$
}
\end{table}
\end{landscape}

\restoregeometry

\end{document}